# Fuzzy models for learning assessment

**Igor Ya. Subbotin[1], Michael Gr. Voskoglou[2]**

[1] Department of Mathematics and Natural Sciences, College of Letters and Sciences, National University, Los Angeles, California, USA
E-mail: isubboti@nu.edu

[2] Department of Mathematical Sciences, School of Technological Applications, Graduate Technological Educational Institute (T. E. I.) of Western Greece, Patras, Greece
E-mail: mvosk@hol.gr

**Abstract**

The concept of learning is fundamental for the study of human cognitive action. In the current paper a Trapezoidal Fuzzy Assessment Model (TFAM) is developed for learning assessment. The TRAMF is a new variation of a special form of the commonly used in Fuzzy Mathematics Center of Gravity (COD) defuzzification technique that we have applied in earlier papers as an assessment method in various human activities. The TFAM's new idea is the replacement of the rectangles appearing in the graph of the COG method by isosceles trapezoids sharing common parts, thus covering the ambiguous cases of students' scores being at the limits between two successive grades (e.g. between A and B). A classroom application is also presented in which the outcomes of the COG and TRAFM methods are compared with those of other traditional assessment methods (calculation of means and GPA index) and explanations are provided for the differences appeared among these outcomes.

**Keywords:** Learning assessment, GPA index, Fuzzy sets, Centre of gravity (COG) defuzzification technique, Trapezoidal Fuzzy Assessment Model (TFAM).

## 1. Introduction

The concept of learning is fundamental for the study of human cognitive action. But while everyone knows empirically what learning is, the understanding of its nature has proved to be complicated. This happens because it is very difficult for someone to understand the way in which the human mind works, and therefore to describe the mechanisms of the acquisition of knowledge by the individual. The problem is getting even harder by taking into consideration the fact that these mechanisms, although they appear to have some common general characteristics, they actually differ in their details from person to person.

There are many theories and models developed by psychologists and education researchers for the description of the mechanisms of learning. Voss [22] argued that learning basically consists of successive problem solving activities, in which the input information is represented of existing knowledge, with the solution occurring when the input is appropriately interpreted.

According to Voss [22] and many other researchers the process of learning involves the following stages: *Representation* of the input data, *interpretation* of this data in order to produce the new knowledge, *generalization* of the new knowledge to a variety of situations and *categorization* of the knowledge. More explicitly the representation of the stimulus input is relied upon the individual's ability to use contents of his/her memory in order to find information that will facilitate a solution development. Learning consists of developing an appropriate number of interpretations and generalizing them to a variety of situations. When the knowledge becomes substantial, much of the process involves categorization, i.e. the input information is interpreted in terms of the classes of the existing knowledge. Thus the individual becomes able to relate the new information to his (her) knowledge structures that have been variously described as schemata, or scripts, or frames.

Voskoglou ([12] and [16 , section 2.3]) developed a stochastic model to describe mathematically the process of learning in the classroom by introducing a finite Markov chain on the stages of the Voss's framework for learning [22] Further, applying principles of the theory of *absorbing Markov chains* [2, Chapter III] on the resulting structure he has obtained a measure for assessing the individuals learning skills [12] and he has also calculated the probability for a student to pass successfully through all the states of the learning process in the classroom [12,15].

However, the knowledge that students have about various concepts is usually imperfect, characterized by a different degree of depth. Also, from the teacher's point of view there exists in many cases vagueness about the degree of his/her students' success in each stage of the learning process. All the above gave us the impulsion in earlier papers to introduce principles of fuzzy logic for a more realistic representation of the process of learning. Thus, Voskoglou presented a fuzzy model for the description of the process of leaning [13] and used the corresponding system's total uncertainty to measure the student's learning skills [14]. Later he also used these ideas in other sectors of mathematical education, like mathematical modelling [18], problem solving [19], etc. On the other hand, Subbotin et al. [4] introduced the idea of applying *the Center of Gravity (COG) defuzzification technique* to learning assessment. Later this idea found interesting continuations and generalizations in the articles [5, 6, 8, 17, 20, 21] , etc. More details about our older researches on fuzzy logic applications are presented in section 2.

Our target in this paper is the expansion of an introduced in [10] *Trapezoidal Fuzzy Model* for learning assessment (TFAM). Accordingly, the rest of the paper is organized as follows: In section 2 we give a brief account of our older research concerning the use of fuzzy logic in the learning process. A particular emphasis is given in this section to the description of a special form of the COG technique, which is actually the basis for the development of the TFAM. In section 3 we describe in detail the TFAM., while in section 4 we present a classroom application illustrating our results in practice. In this application, apart from the fuzzy, we also use traditional methods for learning assessment (calculation of means, GPA index) and we compare their outcomes with those of the COG and the TFAM methods.

For general facts on fuzzy sets and logic we refer to the book [3].

## 2. Fuzzy models for the learning process: Our previous researches

In 1999 Voskoglou developed a fuzzy model for the description of the learning process [13], and later he used the *total uncertainty* of the corresponding fuzzy system for assessing the students' skills in learning mathematics [14]. In Voskoglou's model the major stages of the Voss's framework for learning [22] are represented as fuzzy subsets of a set of linguistic labels characterizing the students' performance and the process of learning is qualitatively studied through the calculation of the possibilities (i.e. their relative membership degrees with respect to the maximum one) of all students' profiles.

Subbotin et al. [4] based on Voskoglou's fuzzy model for learning [13] adapted properly the widely used in Fuzzy Mathematics *Center of Gravity* (COG) *defuzzification technique* (e.g. see [11]) to provide an alternative measure for the learning assessment. Since then, both the authors of the present paper, either collaborating or independently to each other, utilized the COG method for assessing various other students' competencies (e.g. [5], [8], [17], [20], etc), for testing the effectiveness of a CBR system [6] and for assessing the players' performance in the game of Bridge [21].

According to the COG technique the defuzzification of a fuzzy situation's data is succeeded through the calculation of the coordinates of the COG of the level's area contained between the graph of the membership function associated with this situation and the OX axis. In order to be able to design the graph of the membership function we correspond to each $x \in U$ an interval of values from a prefixed numerical distribution, which actually means that we replace $U$ with a set of real intervals. A brief description of the special form of the COG method applied for the learning assessment [4] is the following:

Let $U$ = {A, B, C, D, F} be a set of linguistic labels (grades) characterizing the learners' performance, where A stands for excellent performance, B for very good, C for good, D for mediocre

(satisfactory) and F for unsatisfactory performance respectively. Obviously, the above characterizations are fuzzy depending on the modeler's personal criteria, which however must be compatible to the common logic, in order to model the learning situation in a worthy of credit way. For example, these criteria can be formed by marking the students' performance in the corresponding written or oral test within a scale from 0 to 100 and by assigning to their scores the above characterizations as follows: A = 85-100, B = 84-75, C = 74-60, D =59-50 and F = less than 50. Assume now that we want to assess the general performance of a group, say G, of $n$ students, where $n$ is an integer, $n \geq 2$. For this, we represent G as a fuzzy subset of $U$. In fact, if $n_A$, $n_B$, $n_C$, $n_D$ and $n_F$ denote the number of students that demonstrated excellent, very good, good, mediocre and unsatisfactory performance respectively, we define the *membership function*

$m: U \to [0, 1]$ in terms of the frequencies, i.e. by $m(x) = \frac{n_x}{n}$, for each $x$ in $U$. Then G can be written as a fuzzy subset of $U$ in the form: $G = \{(x, \frac{n_x}{n}): x \in U\}$.

Next we replace $U$ with a set of real intervals as follows: F → [0, 1), D → [1, 2), C → [2, 3), B → [3, 4], A → [4, 5]. Consequently, we have that $y_1 = m(x) = m(F)$ for all $x$ in [0,1), $y_2 = m(x) = m(D)$ for all $x$ in [1,2), $y_3 = m(x) = m(C)$ for all $x$ in [2, 3), $y_4 = m(x) = m(B)$ for all $x$ in [3, 4) and $y_5 = m(x) = m(A)$ for all $x$ in [4,5). Since the membership values of the elements of $U$ in G have been defined in terms of the corresponding frequencies, we obviously have that $\sum_{i=1}^{5} y_i = m(A) + m(B) + m(C) + m(D) + m(F) = 1$ (1).

We are now in position to construct the graph of the membership function $y = m(x)$, which has the form of the bar graph shown in Figure 1. From Figure 1 one can easily observe that the level's area, say S, contained between the bar graph of $y = m(x)$ and the OX axis is equal to the sum of the areas of five rectangles $S_i$, i =1, 2, 3, 4, 5. The one side of each one of these rectangles has length 1 unit and lies on the OX axis.

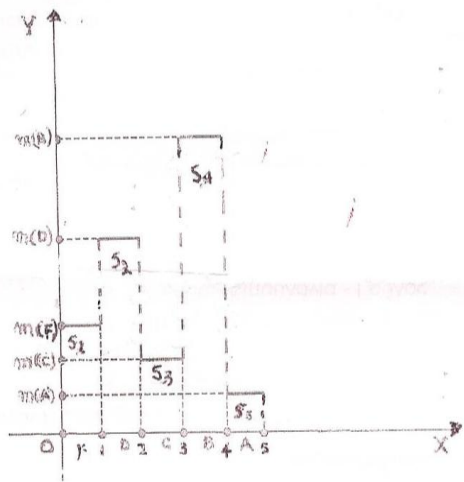

**Figure 1:** Bar graphical data representation

It is well known (e.g. see [23]) that the coordinates $(x_c, y_c)$ of the COG, say $F_c$, of the level's area S are calculated by the formulas:

$$\frac{\iint_S x\, dx dy}{\iint_S dx dy}, \quad \frac{\iint_S y\, dx dy}{\iint_S dx dy} \quad (2).$$

Taking into account the situation presented by Figure 1 and equation (1) it is straightforward to check that in our case formulas (2) can be transformed to the form:

$$x_c = \frac{1}{2}(y_1 + 3y_2 + 5y_3 + 7y_4 + 9y_5),$$

$$y_c = \frac{1}{2}(y_1^2 + y_2^2 + y_3^2 + y_4^2 + y_5^2) \quad (3)$$

In fact, $\iint_S dxdy$ is the area of S which in our case is equal to $\sum_{i=1}^{5} y_i = \frac{n_F + n_D + n_C + n_B + n_A}{n} = 1$

Also $\iint_S xdxdy = \sum_{i=1}^{5}\iint_{S_i} xdxdy = \sum_{i=1}^{5}\int_0^{y_i} dy \int_{i-1}^{i} xdx = \sum_{i=1}^{5} y_i \int_{i-1}^{i} xdx = \sum_{i=1}^{5} y_i [\frac{x^2}{2}]_{i-1}^{i} = \frac{1}{2}\sum_{i=1}^{5} y_i[i^2 - (i-1)^2] = \frac{1}{2}\sum_{i=1}^{5}(2i-1)y_i$

and $\iint_S ydxdy = \sum_{i=1}^{5}\iint_{F_i} ydxdy = \sum_{i=1}^{5}\int_0^{y_i} ydy \int_{i-1}^{i} dx = \sum_{i=1}^{n}\int_0^{y_i} ydy = \frac{1}{2}\sum_{i=1}^{n} y_i^2$

Now, using elementary algebraic inequalities it is easy to check that there is a unique minimum for $y_c$ corresponding to COG $F_m(\frac{5}{2}, \frac{1}{10})$ (an analogous process for TFAM is described in detail in section 3, paragraph 5). Further, the ideal case is when $y_1=y_2=y_3=y_4=0$ and $y_5=1$. Then from formulas (3) we get that $x_c = \frac{9}{2}$ and $y_c = \frac{1}{2}$. Therefore the COG in this case is the point $F_i(\frac{9}{2}, \frac{1}{2})$. On the other hand the worst case is when $y_1=1$ and $y_2=y_3=y_4=y_5=0$. Then from formulas (3) we find that the COG is the point $F_w(\frac{1}{2}, \frac{1}{2})$. Therefore the COG $F_c$ of the level's section F lies in the area of the triangle $F_w F_m F_i$.

Then by elementary geometric observations (the analogous observations are described in detail for TFAM in section 3, paragraphs 5 and 6) one can obtain the following criterion:
- *Between two student groups the group with the bigger $x_c$ performs better.*
- *If the two groups have the same $x_c \geq 2.5$, then the group with the bigger $y_c$ performs better.*
- *If the two groups have the same $x_c < 2.5$, then the group with the lower $y_c$ performs better.*

Recently Subbotin and Bilotskii [7] developed a new variation of the COG method for learning assessment, which they called *Triangular Fuzzy Model* (TFR). The main idea of TFR is to replace the rectangles appearing in the graph of the COG method (Figure 1) by isosceles triangles sharing common parts, so that to cover the ambiguous cases of students scores being at the limits between two successive grades. An improved version of the TFR was applied by Subbotin and Voskoglou [9], for assessing students' critical thinking skills.

### 3. The trapezoidal fuzzy assessment model (TFAM)

The TFAM is a new variation of the presented in the previous section COG method. The novelty of this approach is in the replacement of the rectangles appearing in the graph of the membership function of the COG method (Figure 1) by isosceles trapezoids sharing common parts, so that to cover the ambiguous cases of students scores being at the limits between two successive grades. TFAM was introduced by Subbotin in [10]. Here we shall present an enhanced version of the TFAM. In the TFAM's scheme (Figure 2) we have five trapezoids, corresponding to the students' grades F, D, C, B and A respectively defined in the previous section. Without loss of generality and for making our calculations easier we consider isosceles trapezoids with bases of length 10 units lying on the OX axis. The height of each trapezoid is equal to the percentage of students who achieved the corresponding grade, while the parallel to its base side is equal to 4 units. We allow for any two adjacent trapezoids to have 30% of their bases (3 units) belonging to both of them. In this way we cover the ambiguous cases of students' scores being at the limits between two successive grades. It is a very common approach to divide the interval of the specific grades in three parts and to assign the corresponding grade using + and - . For example, 75 – 77 = B-, 78 – 81 = B, 82 – 84 = B+. However, this consideration does not reflect the common situation, where the teacher is not sure about the grading of the students whose performance could be assessed as marginal between and close to two adjacent grades; for example, something like 84 - 85 being between $B_+$ and $A_-$. The TFAM fits this situation.

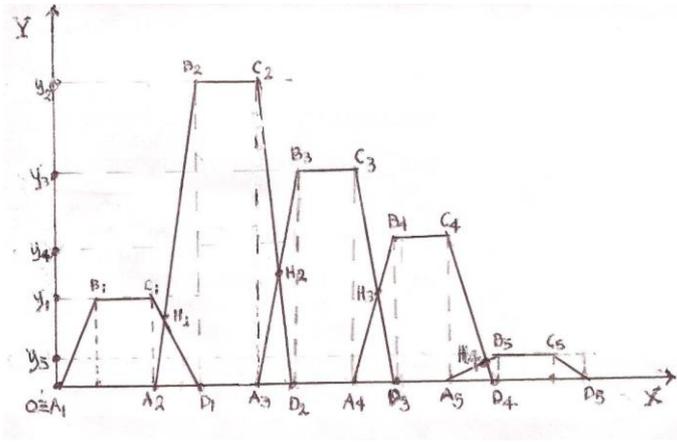

**Figure 2:** The TRAFM's scheme

A student group can be represented, as in the COG method, as a fuzzy set in *U*, whose membership function y=m(x) has as graph the line $OB_1C_1H_1B_2C_2H_2B_3C_3H_3B_4C_4H_4B_5C_5D_5$ of Figure 2, which is the union of the line segments $OB_1$, $B_1C_1$, $C_1H_1$,……., $B_5C_5$, $C_5D_5$. However, in case of the TRAFM the analytic form of y = m(x) is not needed for calculating the COG of the resulting area. In fact, since the marginal cases of the students scores are considered as common parts for any pair of the adjacent trapezoids, it is logical to count these parts twice; e.g. placing the ambiguous cases B+ and A- in both regions B and A. In other words, the COG method, which calculates the coordinates of the COG of the area between the graph of the membership function and the OX axis, thus considering the areas of the "common" triangles $A_2H_1D_1$, $A_3H_2D_2$, $A_4H_3D_3$ and $A_5H_4D_4$ only once, is not the proper one to be applied in the above situation.

Instead, in this case we represent each one of the five trapezoids of Figure 2 by its COG $F_i$, i=1, 2, 3, 4, 5 and we consider the entire area, i.e. the sum of the areas of the five trapezoids, as the system of these points-centers. More explicitly, the steps of the whole construction of the TRAFM are the following:

    1. Let $y_i$, i=1, 2, 3, 4, 5 be the percentages of the students whose performance was characterized by the grades E, D, C, B, and A respectively; then $\sum_{i=1}^{5} y_i = 1$ (100%).

    2. We consider the isosceles trapezoids with heights being equal to $y_i$, i=1, 2, 3, 4, 5, in the way that has been illustrated in Figure 2.

    3. We calculate the coordinates ($x_{c_i}$, $y_{c_i}$) of the COG $F_i$, i=1, 2, 3, 4, 5, of each trapezoid as follows: It is well known that the COG of a trapezoid lies along the line segment joining the midpoints of its parallel sides *a* and *b* at a distance *d* from the longer side *b* given by $d = \frac{h(2a+b)}{3(a+b)}$, where *h* is its height (e.g. see [24])..Therefore in our case we have $y_{c_i} = \frac{y_i(2*4+10)}{3*(4+10)} = \frac{3y_i}{7}$. Also, since the abscissa of the COG of each trapezoid is equal to the abscissa of the midpoint of its base, it is easy to observe that $x_{ci}=7i-2$.

    4. We consider the system of the COG's $F_i$, i=1, 2, 3, 4, 5 and we calculate the coordinates ($X_c$, $Y_c$) of the COG $F_c$ of the whole area *S* considered in Figure 2 by the following formulas, derived from the commonly used in such cases definition (e.g. see [25]): $X_c = \frac{1}{S}\sum_{i=1}^{5} S_i x_{c_i}$, $Y_c = \frac{1}{S}\sum_{i=1}^{5} S_i y_{c_i}$ (4).

In formulas (4) $S_i$, i= 1, 2, 3, 4, 5 denote the areas of the corresponding trapezoids. Thus, $S_i = \frac{(4+10)y_i}{2} = 7y_i$ and $S = \sum_{i=1}^{5} S_i = 7\sum_{i=1}^{5} y_i = 7$. Therefore, from formulas (4) we finally get that

$$X_c = \frac{1}{7}\sum_{i=1}^{5} 7y_i(7i-2) = (7\sum_{i=1}^{5} iy_i) - 2$$

(5)

$$Y_c = \frac{1}{7}\sum_{i=1}^{5} 7y_i(\frac{3}{7}y_i) = \frac{3}{7}\sum_{i=1}^{5} y_i^2$$

5. We determine the area where the COG $F_c$ lies as follows: For $i, j=1, 2, 3, 4, 5$, we have that $0 \leq (y_i - y_j)^2 = y_i^2 + y_j^2 - 2y_iy_j$, therefore $y_i^2 + y_j^2 \geq 2y_iy_j$, with the equality holding if, and only if, $y_i = y_j$. Therefore

$$1 = (\sum_{i=1}^{5} y_i)^2 = \sum_{i=1}^{5} y_i^2 + 2\sum_{\substack{i,j=1,\\ i\neq j}}^{5} y_iy_j \leq \sum_{i=1}^{5} y_i^2 + 2\sum_{\substack{i,j=1,\\ i\neq j}}^{5} (y_i^2 + y_j^2) = 5\sum_{i=1}^{5} y_i^2 \text{ or } \sum_{i=1}^{5} y_i^2 \geq \frac{1}{5} \quad (6),$$

with the equality holding if, and only if, $y_1 = y_2 = y_3 = y_4 = y_5 = \frac{1}{5}$. In the case of equality the first of formulas (5) gives that $X_c = 7(\frac{1}{5} + \frac{2}{5} + \frac{3}{5} + \frac{4}{5} + \frac{5}{5}) - 2 = 15$. Further, combining the inequality (6) with the second of formulas (5) one finds that $Y_c \geq \frac{3}{35}$. Therefore the unique minimum for $Y_c$ corresponds to the COG $F_m(15, \frac{3}{35})$. The ideal case is when $y_1 = y_2 = y_3 = y_4 = 0$ and $y_5 = 1$. Then from formulas (5) we get that $X_c = 33$ and $Y_c = \frac{3}{7}$. Therefore the COG in this case is the point $F_i(33, \frac{3}{7})$. On the other hand, the worst case is when $y_1 = 1$ and $y_2 = y_3 = y_4 = y_5 = 0$. Then from formulas (3), we find that the COG is the point $F_w(5, \frac{3}{7})$. Therefore the area where the COG $F_c$ lies is the area of the triangle $F_w F_m F_i$ (see Figure 3).

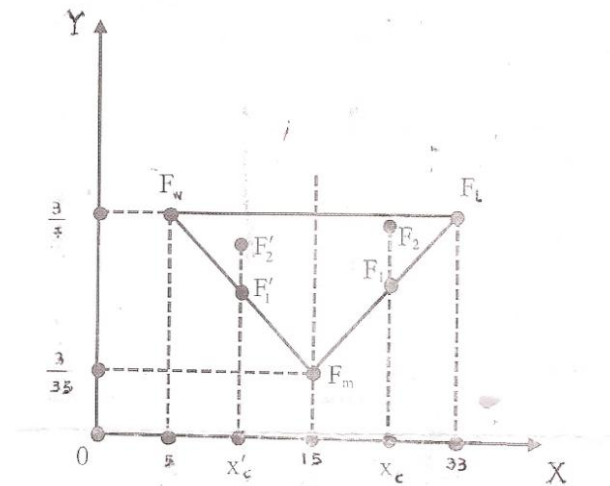

**Figure 3:** The area where the COG lies

6. We formulate our criterion for comparing the performances of two (or more) different student groups' as follows: From elementary geometric observations (see Figure 3) it follows that for two student groups the group having the greater $X_c$ performs better. Further, if the two groups have the same $X_c \geq 15$, then the group having the COG which is situated closer to $Fi$ is the group with the greater $Y_c$. Also, if the two groups have the same $X_c < 15$, then the group having the COG which is situated farther to $Fw$ is the group with the smaller $Y_c$. Based on the above considerations it is logical to formulate our criterion for comparing the two groups' performance in the following form:

- *Between two student groups the group with the greater value of $X_c$ demonstrates the better performance.*
- *If two student groups have the same $X_c \geq 15$, then the group with the greater value of $Y_c$ demonstrates the better performance.*
- *If two student groups have the same $X_c < 15$, then the group with the smaller value of $Y_c$ demonstrates the better performance.*

**4. A classroom application**

The students of two different Departments of the School of Technological Applications (prospective engineers) of the Graduate T. E. I. of Western Greece acheived the following scores in the final common exam of the mathematics course of their first term of studies (the contents of the course were the same and the instructor was the same person for both Departments):
*Department 1* ($D_1$)*:* 99(1 student), 83(2), 82(1), 74(10), 72(2), 70(1), 59(10), 55(2), 48(7), 45(2).
*Department 2* ($D_2$): 85(2), 75(1), 62(2), 60(10), 52(1), 50(8), 25(4), 10(1).
The results of the students' performance are summarized in Table 1 below:

**Table 1**: The exam's results

| Grade | $D_1$ | $D_2$ |
|---|---|---|
| A | 1 | 2 |
| B | 3 | 1 |
| C | 13 | 12 |
| D | 12 | 9 |
| F | 9 | 5 |
| Total | 38 | 29 |

The evaluation of the above results will be performed below using both traditional methods, based on principles of the classical (bivalent) logic, and fuzzy logic methods.

### 4.1 Traditional methods

*i) Calculation of the means:* A straightforward calculation gives that the means of the students' scores are approximately 62.231 and 52.793 for $D_1$ and $D_2$ respectively. This shows that the *mean performance* was good (C) for the students of $D_1$ and satisfactory (D) for the students of $D_2$.

*ii) Calculation of the GPA index:* We recall that the *Great Point Average* (*GPA*) index is a weighted mean, where more importance is given to the higher scores by attaching greater coefficients (weights) to them (e.g. see [1]). In other words, the GPA index focuses on the *quality performance* of a student group.

Let us denote by $n_A, n_B, n_C, n_D$ and $n_E$ the numbers of students of a given group whose performance is characterized by A, B, C, D and F respectively and by $n$ the total number of students of the group..

Then, the GPA index is calculated by the formula GPA=$\frac{0n_F + n_D + 2n_C + 3n_B + 4n_A}{n}$. Obviously we have that $0 \leq$ GPA $\leq 4$.

In our case the above formula can be written as GPA $= y_2 + 2y_3 + 3y_4 + 4y_5$ (7). Then, using the data of Table 1 it is easy to check that the GPA for $D_1$ is equal to $\frac{51}{38} \approx 1.342$ and for $D_2$ is equal to $\frac{44}{29} \approx 1.517$. Therefore, since the values of the GPA index are less than the half of its maximal possible value, which is equal to 4, the quality performance of both Departments was less than satisfactory. However, in contrast to their mean performances, the quality performance of $D_2$ was better than the corresponding performance of $D_1$.

### 4.2 Fuzzy logic methods

In this paragraph we shall apply the fuzzy logic methods described in sections 2 and 3 of this paper as follows:

*iii) The COG method:* Observing the coefficients of the $y_i$'s, i=1, 2, 3, 4, 5, in the first of formulas (3) and taking into account that, according to the criterion stated in section 2, the COG's abscissa $x_c$ measures a student group's performance, it becomes evident that the COG method is also focused, as the GPA index does, on the student groups' *quality performance*.

In case of our classroom application taking into account the data of Table 1 and using the first of formulas (3) we find that $x_c = \frac{140}{76} \approx 1.842$ for $D_1$ and

$x_c = \frac{117}{58} \approx 2.017$ for $D_2$. Since the above values of $x_c$ are less than the half of its value in the ideal case, which is equal to $\frac{9}{2}$ (see section 2), the quality performance of both Departments according was less than satisfactory. Further $D_2$ demonstrated a better quality performance than $D_1$.

*iv) Application of TFAM:* Observing the coefficients of the $y_i$'s, i=1, 2, 3, 4, 5, in the first of formulas (5) it is easy to conclude that the COG method is also focused on the student groups' *quality performance*.

In case of our classroom application taking into account the data of Table 1 and using the first of formulas (5) we find that $X_c = \frac{623}{39} \approx 15.974$ for $D_1$ and $X_c = \frac{511}{29} \approx 17.621$ for $D_2$. In this case the value of $X_c$ for $D_1$ is less than the half of its value in the ideal case, which is equal to 33 (see section 3). This shows that the quality performance of $D_1$ according was less than satisfactory. On the contrary, since the value of $X_c$ for $D_2$ is greater than the half of its value in the ideal case, the quality performance of $D_2$ was more than satisfactory.

**4.3 Comparison of the assessment methods used**

In paragraphs 4.1 and 4.2 we have applied four in total methods for learning assessment. The first of these methods measures the mean performance of a student group, while the other three methods (GPA, COG and TFAM) measure its quality performance by assigning greater coefficients (weights) to the higher scores. The coefficients attached to the $y_i$'s in these three methods -see formula (7) and the first of formulas (3) and (5) respectively- are present in the following Table 2:

**Table 2:** Weight coefficients of the $y_i$'s

| $y_i$ | GPA | COG ($x_c$) | TRAFM ($X_c$) |
|---|---|---|---|
| $y_1$ | 0 | 1/2 | 7 |
| $y_2$ | 1 | 3/2 | 14 |
| $y_3$ | 2 | 5/2 | 21 |
| $y_4$ | 3 | 7/2 | 28 |
| $y_5$ | 4 | 9/2 | 35 |

From Table 2 becomes evident that TFAM assigns greater coefficients to the higher with respect to the lower scores than COG and also COG does the same thing with respect to GPA. In other words TFAM is more accurate than COG, and COG is more accurate than GPA for measuring the quality performance of a student group.

This explains why in our classroom application the quality performances of $D_1$ and $D_2$ were found to be less than satisfactory by applying the GPA and COG methods, while the application of the TRAFM method has demonstrated a more than satisfactory quality performance for $D_2$ and less than satisfactory for $D_1$.

One should also mention that, while $D_2$ demonstrated in all cases (GPA, GOG and TRAFM) a better quality performance than $D_1$, in contrast to the mean performance of $D_1$, which was found to be better than the corresponding performance of $D_2$ (first method of paragraph 4.1).

In concluding, it is suggested to the user of the above four assessment methods to choose the one that fits better to his/her personal goals.

**5. Conclusion and discussion**

The methods for assessing a group's performance (for human activities) usually applied in practice are based on principles of the bivalent logic (yes-no). However, fuzzy logic, due to its nature of characterizing a situation with multiple values by using linguistic variables, offers a wider and richer field of resources for this purpose. This gave us the impulsion to introduce in this paper principles of fuzzy logic for developing an expanded version of the TRAFM approach for learning assessment. The TRAFM is actually a more sensitive version of the COG fitting better to the ambiguous cases of students' scores lying at the limits between two different grades. We also presented a classroom application in which we have compared the outcomes of TRAFM approach with the corresponding outcomes of the COG technique and of other traditional assessment methods (calculation of the means and GPA index).

However, there is a need for more classroom applications to be performed in future for obtaining safer statistical data. On the other hand, since the TRAFM approach appears to have the potential of a general assessment method, our future research plans include also the effort to apply this approach for assessing the individuals' performance in several other human activities.